\documentclass[reqno, 12pt]{amsart}

\newtheorem{theorem}{Theorem}

\newtheorem{corollary}[theorem]{Corollary}
\newtheorem{definition}{Definition}

\newtheorem*{problem}{Problem}

\usepackage{amssymb}
\usepackage{amsmath}
\usepackage{amsfonts}

\newcommand{\degr}{{\text{deg}}}
\newcommand{\Cantor}{{\mathcal C}}

\begin{document}

\title{On homeomorphic product measures on the Cantor set}

\author{Randall Dougherty and R. Daniel Mauldin}
\address{IDA Center for Communications Research, 4320 Westerra Ct.,
San Diego, CA 92121}
\email{rdough@ccrwest.org}
\address{Department of Mathematics, PO Box 311430, University of North Texas,
Denton, TX 76203}
\email{mauldin@unt.edu}


\begin{abstract}

Let $\mu(r)$ be the Bernoulli measure on the Cantor space given as
the infinite product of two-point measures with weights
$r$ and $1-r$.  It is a long-standing open problem
to characterize those $r$ and $s$ such that $\mu(r)$
and $\mu(s)$ are topologically equivalent (i.e., there
is a homeomorphism from the Cantor space to itself sending
$\mu(r)$ to $\mu(s)$).  The (possibly) weaker property of $\mu(r)$
and $\mu(s)$ being continuously reducible to each other
is equivalent to a property of $r$ and $s$ called
binomial equivalence.
In this paper we define an algebraic property called
``refinability'' and show that, if $r$ and $s$ are
refinable and binomially equivalent, then $\mu(r)$ and
$\mu(s)$ are topologically equivalent.
We then give a class of examples of refinable numbers;
in particular, the positive numbers $r$ and $s$ such
that $s = r^2$ and $r = 1-s^2$ are refinable, so the
corresponding measures are topologically equivalent.

\end{abstract}

\thanks{The second author was supported in part by NSF grant: DMS 0100078
}
\maketitle


Two measures $\mu$ and $\nu$ defined on the family of Borel subsets of
a topological space $X$ are said to be {\it homeomorphic} or
{\it topologically equivalent} provided there exists a homeomorphism
$h$ of $X$ onto
$X$ such that $\mu$ is the image measure of $\nu$ under $h$:
$\mu = \nu h^{-1}.$ This means $\mu(E) = \nu(h^{-1}(E))$ for each
Borel subset $E$ of $X$.

One may be interested in the structure of
these equivalence classes of measures or in a particular equivalence
class. For example, a probability measure $\mu$ on $[0,1]$ is
topologically equivalent to Lebesgue measure if and only if $\mu$ gives
every point measure 0 and every non-empty open set positive measure.
(The distribution function of $\mu$ is a
homeomorphism on $[0,1]$ witnessing this equivalence.)
This is a special case of a result of
Oxtoby and Ulam \cite{OU}, who characterized those probability
measures $\mu$ on finite dimensional cubes $[0,1]^n$ which are
homeomorphic to Lebesgue measure. For this to be so, $\mu$
must give points measure 0, non-empty open sets positive measure, and
the boundary of the cube measure 0. Later Oxtoby and Prasad \cite{OP} extended
this result to the Hilbert cube. These results have been
extended and applied to various manifolds. The book of Alpern and
Prasad \cite{AP} is an excellent source for these developments. Oxtoby 
\cite{O} also
characterized those measures on the space of irrational numbers in $[0,1]$ 
which
are homeomorphic to Lebesgue measure.

It is natural to ask
what measures are homeomorphic to Lebesgue measure on $\Cantor=\{0,1\}^N$,
   the Cantor space, where by Lebesgue measure we mean Haar measure or
   infinite product measure $\mu({1/2})$ resulting from fair coin
   tossing.  The topology on $\Cantor$ is the standard product topology;
   we will use as basic open (actually clopen) sets for this topology
   the sets $\langle e \rangle$ for all finite sequences $e$ from $\{0,1\}$,
   where $\langle e \rangle$ is the set of infinite sequences in $\Cantor$
   which begin with the finite sequence $e$.  (These basic clopen
   sets are sometimes called cylinders.)  We will say that the
   {\it length} of a basic clopen set $\langle e \rangle$ is the
   length of the finite sequence $e$.

   It turns out that the Cantor space is more rigid than
   $[0,1]^n$ for measure homeomorphisms -- it is not true
   that a measure $\nu$ on $\Cantor$ which gives points measure $0$ and
   non-empty open sets positive measure is equivalent to Lebesgue
   measure. In fact, even among the product measures the only one which
   is equivalent to Lebesgue measure is Lebesgue measure itself. To
   describe the situation let us use the following notation. For each
   number $r$, $0 \leq r \leq 1$, let $\mu(r)$ be the infinite product
   measure determined by coin tossing with probability of success $r$.
Consider the equivalence relation on $[0,1]$,  $r \sim_{top} s$ if and
only if $\mu(r)$ is topologically equivalent to $\mu(s)$.  (We will
sometimes abuse terminology by saying that $r$ is topologically
equivalent to $s$.)
It turns
out that this equivalence relation is closely related to an
algebraic/combinatorial
relation. To explain this, we make the following definition.



\begin{definition}
\label{binreddef}
Let  \ $0 < r,s < 1$. The number $s$ is said to be
     binomially reducible to $r$ provided
\begin{equation}\label{brfmla}
s = \sum_{i=0}^n \ a_i \ r^i(1-r)^{n-i},
\end{equation}
where $n$ is a non-negative integer and each $a_i$ is an integer with $0
\leq a_i \leq \binom {n}{i}$.
\end{definition}

The numbers $r^a(1-r)^b$ for integers $a,b\ge 0$ will be referred
to as {\it cylinder sizes} for $r$; they are the measures of the
basic clopen sets under $\mu(r)$.  So the right side of \ref{brfmla}
is the general form for the measure of a clopen set under $\mu(r)$.

Let us note some basic facts about reducibility to be used later.
(See \cite{Ma} for this and further background information.)
Note that if $s$ is
binomially reducible to $r$, so is $1-s$ (change $a_i$ to $\binom
{n}{i} - a_i$). If
$s_1$ and $s_2$ are reducible to $r$, so is $s_1s_2$.
(If $s_1 =
   \sum_{i=0}^n a_ir^i(1-r)^{n-i}$ and $s_2 =
   \sum_{j=0}^m b_jr^j(1-r)^{m-j}$, then $s_1s_2 =
   \sum_{i=0}^n \sum_{j=0}^m a_ib_jr^{i+j}(1-r)^{n+m-i-j}
   =\sum_{k=0}^{n+m}
  (\sum_{i+j=k} a_ib_j)r^{k}(1-r)^{n+m-k}$ and $ \sum_{i+j=k} a_ib_j
   \leq \sum_{i+j=k} \binom {n}{i}\binom{m}{j} = \binom {n+m}{k}$.)
Hence, if $s$ is binomially reducible to $r$, so is $s^a(1-s)^b$ for
any $a,b \geq 0$. Also, it is known that $\mu(s)$ is continously reducible
to or is a continuous
image of the measure $\mu(r)$ (i.e., $\mu(s) = \mu(r) \circ g^{-1}$ for
some continuous $g: \Cantor \to \Cantor$) if and only if $s$ is binomially
reducible to $r$ \cite{Ma}. Thus, we have another natural equivalence relation
on $[0,1]$.

\begin{definition}
\label{bineqdef}
Let  \ $0 < r,s < 1$. Then $r$ is binomially
   equivalent to $s$, denoted $r \approx s$, provided $r$ is binomially
     reducible to $s$ and $s$ is binomially reducible to $r$, or,
     equivalently, each of the measures $\mu(r)$ and $\mu(s)$ is a
     continuous image of the other.
\end{definition}

Among several still unsolved problems concerning these relations
is the following.

\begin{problem}[{\cite[Problem 1065]{Ma}}]
Is it true that the product measures $\mu(r)$ and $\mu(s)$ are homeomorphic if
and only if each is a continuous image of the other, or, equivalently,
each of the numbers $r$ and $s$ is binomially reducible to the other?
\end{problem}

(Note: After this paper was first circulated, the above problem
was solved in the negative by Austin~\cite{Au}.)

One can think of this problem in the following way. Suppose we have
$\mu(s) = \mu(r)\circ g^{-1}$ and $\mu(r) = \mu(s) \circ h^{-1}$, where 
the maps $g$
and $h$ are continuous. Is there some
sort of Cantor-Bernstein or back-and-forth argument for the
Cantor set  which, given $g$ and $h$,
produces not just a one-to-one onto map, but a
homeomorphism taking $\mu(r)$ to $\mu(s)$?

Many cases of this problem have already been settled.
Let us note that, for a given $n$,
the functions $r^i(1-r)^{n-i}$ for $0 \leq i \leq n$
are linearly independent polynomials, since their trailing terms
(i.e., nonzero terms of least degree) have distinct degrees. Therefore,
$\sum_{i=0}^na_ir^i(1-r)^{n-i}$ as in Definition~\ref{binreddef}
is a polynomial of degree $> 1$ unless
it is $0$ (when $a_i = 0$ for all $i$), $1$ (when $a_i = \binom {n}{i}$),
$1-r$ ($a_i =
\binom {n-1}{i}$),
or $r$ ($a_i = \binom {n-1} {i-1}$). Therefore, if $r$ and $s$ are binomially
   reducible to each other, $r\neq s$, and $r \neq 1-s$, then $s = P(r)$
   and $r = Q(s)$ where $\degr(P), \degr(Q)>1$, so $r = Q\circ P(r)$ and
   $\degr(Q\circ P)
   > 1$. Thus, $r$ is algebraic. Also, in this case, $r$ and $s$ have
   the same algebraic degree. Moreover, $r$ is an algebraic integer
if and only if $s$ is. Huang~\cite{H} showed that if $r$ is an algebraic 
integer of
degree 2, and $r \approx s$, then $r = s$ or $r = 1-s$.
In fact, Navarro-Bermudez~\cite{NB} showed that if $r$ is rational or
transcendental and $r\approx s$, then $r =s$ or $r = 1-s$.
We gather these facts in the following theorem.

\begin{theorem}[various authors] For $r$ rational, transcendental,
or an algebraic integer of degree
$2$, the $\sim_{top}$ equivalence class containing $r$ and the $\approx$
equivalence class containing $r$ are both equal to $\{r,1-r\}.$
\end{theorem}

On the other hand, it is known that for every $n \geq 3$, there are
algebraic integers $r$ of degree $n$ such that the $\approx$
equivalence class containing $r$ has at least 4 elements \cite
{H}. (In fact, Pinch \cite{P} showed that, if $n = 2^{k+1}$, then there is an
algebraic integer $r$ of degree $n$ with at least $2k$ distinct numbers
binomially equivalent to it.) The simplest
of these is the solution of
  $$
r^3 + r^2 -1 =0
$$
lying in the open interval $(0,1)$. For this value of $r$, it turns out
that $s = r^2 \approx r$, and Navarro-Bermudez and Oxtoby \cite {ONB} 
proved that
$r\sim_{top} s$ via a simple homeomorphism. Until now this has been the 
only nontrivial example of
topologically equivalent product measures.

The purpose of this paper is to present a new condition under which
binomially equivalent numbers are topologically equivalent. First, we 
define a condition
called ``refinable'' on  numbers in $[0,1]$. Next, we show that  if $r$ 
and $s$ are
binomially equivalent and both $r$ and $s$ are refinable, then the
measures $\mu(r)$ and $\mu(s)$ are homeomorphic. Finally, we apply our
condition to the root $r$ of
$$
r^4 +r -1 = 0,
$$
with $r$ between 0 and 1, and to $s = r^2$. We show both $r$ and $s$
are refinable
and $r$ and $s$ are binomially equivalent. Thus, $\mu(r)$ and $\mu(s)$
are topologically equivalent via a very non-trivial homeomorphism.

Any cylinder size $r^a(1-r)^b$ can be split into two cylinder sizes
$r^{a+1}(1-r)^b$ and $r^a(1-r)^{b+1}$.  Either or both of these
can be split in the same way, and so on.  After finitely many
steps, one has partitioned the original cylinder size into finitely
many cylinder sizes.  We will call a partition obtained in this
way a {\it tree partition} of $r^a(1-r)^b$ (named after the
representation of the Cantor space as the set of paths through a
complete infinite binary tree).  A tree partition corresponds
to a partition of a basic clopen set in $\Cantor$ into basic clopen subsets.
Note that any tree partition
can be split further by steps as above to yield a new tree
partition in which
all the final cylinder sizes have
the same length, say $a+b+n$; in this final partition the
cylinder size $r^{a+i}(1-r)^{b+n-i}$ will occur $\binom ni$ times.

On the other hand,
one may be able to partition a cylinder size for $r$
into a finite collection of smaller cylinder sizes (whose sum is the
original cylinder size; repetitions are allowed)
in a way which is not a tree partition.  For instance, one can partition
the cylinder size $1$ into $\{r^3,(1-r)^3,r(1-r),r(1-r),r(1-r)\}$
(or, written more briefly, $\{r^3,(1-r)^3,3r(1-r)\}$; we will treat
a positive integer coefficient of a cylinder size as a multiplicity).
For specific values of $r$, there may be many more such partitions.

Recall the definition of refinement:
given partitions $P$ and $P'$ of the same set, we say that $P'$ is
a {\it refinement} of $P$ (or $P$ has been {\it refined} to $P'$)
if every member of $P$ is a union of members of $P'$.  The corresponding
definition for partitions of a number (e.g., a cylinder size) is:
$P'$ is a refinement of $P$ if one can write $P'$ as the union
(respecting multiplicities) of collections $S_t$ for $t \in P$ such
that, for each $t$ in $P$, the sum of $S_t$ is $t$.

\begin{definition}
A number $r$ is refinable provided every partition of a cylinder size
for $r$ into smaller cylinder sizes
can be refined to a tree partition.\end{definition}

An equivalent definition in symbols: $r$ is refinable iff, for every
true equation of the form
$$
r^a(1-r)^b = r^{c_1}(1-r)^{d_1} + r^{c_2}(1-r)^{d_2}+ \ldots
+r^{c_m}(1-r)^{d_m},
$$
there exist $n \geq 0$ and nonnegative integers $p_{ij}$ for $0 \leq
i \leq n$, $1 \leq j \leq m$ such that
$$
\binom {n}{i} = p_{i1}+p_{i2}+\ldots+p_{im}
$$
for $0\leq i \leq n$ and
$$
r^{c_j}(1-r)^{d_j} = \sum_{i=0}^n p_{ij}r^{a+i}(1-r)^{b+n-i}
$$
for $1 \leq j \leq m.$
We note that in this definition $a$, $b$, $c_j$, and $d_j$ are assumed to be
nonnegative integers, but the definition would be equivalent if we
allowed them to be arbitrary integers.

We briefly compare this
notion to that of a ``good'' measure as introduced by Akin \cite{Ak}.
A probability measure $\mu$ on the Cantor space is good if,
whenever $U,V$ are clopen sets with $\mu(U) < \mu(V)$, there
exists a clopen subset $W$ of $V$ such that $\mu(W) = \mu(U)$.
We state a few facts without proof here.
If a product measure $\mu(r)$ is good, then $r$ is refinable.
If $r$ is transcendental, then $r$ is refinable, but $\mu(r)$ is
not good. If $r$ is rational and $r\neq 1/2$, then $r$ is not
refinable and hence $\mu(r)$ is not good.

Refinability is useful because of the following result.

\begin{theorem}
\label{newthm}
If $0 < r,s < 1$, $r$ and  $s$ are binomially equivalent, and each of $r$ 
and $s$ is refinable,
then the measures $\mu(r)$ and $\mu(s)$ are homeomorphic.
\end{theorem}

\begin{proof} We construct partitions $P_n$ and $Q_n$ of $\Cantor$ into
   clopen sets for $n = 0,1,2,\ldots$ and bijections $\pi_n:P_n\mapsto
   Q_n$ satisfying the following properties:
\begin{enumerate}
\item $P_{n+1}$ is a
   refinement of $P_n$ and $Q_{n+1}$ is a refinement of $Q_n$,
\item each member of $P_{2n-1}$ and each member of $Q_{2n}$ is a basic
   clopen set of length $\geq n$,
\item for any $X \in P_n$ we have
   $\mu(s)(\pi_n(X)) = \mu(r)(X)$, and
\item if $X\in P_{n+1}$ and $X
   \subseteq X' \in P_n$, then $\pi_{n+1}(X) \subseteq \pi_n(X')$.
\end{enumerate}

Given the above sequence, define $f:\Cantor\mapsto \Cantor$ by: for each
$\alpha \in \Cantor$, let $X_n$ be the unique member of $P_n$ containing
$\alpha$ and let $f(\alpha)$ be the unique element of $\bigcap_n
\pi_n(X_n).$ It is straightforward to verify that $f$ is a well-defined
homeomorphism of $\Cantor$ ($f^{-1}$ is defined by an analogous method from
$Q_n$ to $P_n$), and $f(X) = \pi_n(X)$ for all $X\in P_n,$ so that
$\mu(s)(f(X))=\mu(r)(X)$ for $X\in \bigcup_nP_n$. Since every clopen set
$A$ is a finite disjoint union of sets each in $\bigcup_n P_n$, $f$ maps
$\mu(r)$ to $\mu(s)$.

We build $P_n$, $Q_n$, and $\pi_n$ by a back-and-forth recursive
construction.  Let
$P_0 = Q_0 = \{\Cantor\}$ with $\pi_0(\Cantor) = \Cantor$. Given
$P_{2n},Q_{2n},\pi_{2n},$ let $P_{2n+1}$ be a refinement of $P_{2n}$
into basic clopen sets of length $\geq n+1$. Fix $Y\in Q_{2n}$, a basic
clopen set, say of $\mu(s)$-measure $s^a(1-s)^b$. Now,
$\pi^{-1}_{2n}(Y)\in P_{2n}$ is a union of basic clopen sets
$X_1,...,X_k \in P_{2n+1}$, each having $\mu(r)$-measure $r^p(1-r)^q$
for some integers $p,q$, and these measures add up to $s^a(1-s)^b$.
Since $r$ is binomially reducible to $s$, so is each $r^p(1-r)^q$. Thus,
each $\mu(r)(X_j)$ can be expressed as a finite sum of numbers
$s^c(1-s)^d$. Putting these together for all such $X_j$'s, we get a
list of numbers $s^{c_1}(1-s)^{d_1}, s^{c_2}(1-s)^{d_2},\ldots,
s^{c_m}(1-s)^{d_m}$ with sum $s^a(1-s)^b$. Since $s$ is refinable, $Y$
can be partitioned into clopen sets $\hat{Y}_1,\ldots,\hat{Y}_m$ with
$\mu(\hat{Y}_i) =s^{c_i}(1-s)^{d_i}$ for each $i$. But the list above
was obtained by
joining the lists for the individual $X_j$'s together; hence, we can
combine the $\hat{Y}_i$'s to get clopen sets $Y_1, Y_2,\ldots,Y_k$, still 
forming a
partition of $Y$, such that $\mu(s)(Y_j) = \mu(r)(X_j)$, for each $j$.
Put these sets $Y_j$ into $Q_{2n+1}$, letting $\pi_{2n+1}(X_j) = Y_j$.
Once this is done for all $Y\in Q_{2n}$, we will have the desired
partition $Q_{2n+1}$ and map $\pi_{2n+1}$.

We have finished refining the partition on the $P$ side; it is now
the partition on the $Q$ side that needs to be refined next.
So let $Q_{2n+2}$ be a refinement of $Q_{2n+1}$ into basic clopen sets
of length $\geq n+1$, and apply the above procedure with $r$ and
$s$ interchanged to get $P_{2n+2}$ and $\pi_{2n+2}$ (the map from
$Q_{2n+2}$ to $P_{2n+2}$ will be $\pi_{2n+2}^{-1}$).  This will
complete the back-and-forth recursive step.
\end{proof}

To prove that a number $r$ is refinable, it suffices to show the following.
Given any two finite multisets (sets whose elements can have
multiplicity greater than $1$) $A$ and $B$ of cylinder sizes for $r$ such
that $\sum A = \sum B$, one can transform $A$ and $B$ to a common
multiset $C$, where for $B$ the allowed transform steps are arbitrary
splits (meaning replace a cylinder size $x$ with any collection of
cylinder sizes that add to $x$), while for $A$ the only allowed steps
are tree splits (meaning replace $x$ with $xr$ and $x(1-r)$). If
$A$ is a singleton, $A = \{r^a(1-r)^b\}$, and $B$ and $C$ are as just
described, then $C$ will be a tree partition of $r^a(1-r)^b$ which is
a refinement of the partition $B$.
[Note that we only need the case where $A$ is a singleton, but
proofs that one can transform $A$ and $B$ to $C$ as above will usually
work even when $A$ is an arbitrary multiset.]

Equivalently, one can
transform $A$ and $B$ to a common $C'$ where one gets from $B$ to $C'$
by arbitrary splits and one get from $A$ to $C'$ by a sequence of tree
splits followed by a sequence of merges (meaning replace a subcollection
of the current multiset by its sum, which we may or may not require to
be a cylinder size). This is because, if $C$ is the multiset obtained from
$A$ by the tree splits alone, then $C$ is obtained from $C'$ by
arbitrary splits
(a split is the inverse of a merge) and hence from $B$ by arbitrary splits.

In fact, it will suffice to transform $A$ to $C'$ by any sequence of
merges and tree splits in any order, because a merge followed by a
tree split is equivalent to one or more tree splits followed by a merge
(if $\sum_ix_i = x$, then $\sum_ix_ir = xr$ and
$\sum_ix_i(1-r) = x(1-r)$). Define a ``tree move'' to be a merge or a
tree split.
This is the method we will use in the refinability proofs to follow:
given $A$ and $B$, transform $A$ to $A'$ by tree moves and
$B$ to $B'$ by splits, and show that $A' = B'$ (this is the common
multiset $C'$).

We will describe such transformations as built up out of simple steps.
For instance, suppose we have cylinder sizes $a_1$, $a_2$, $b_1$, $b_2$,
and $b_3$, and we demonstrate how to transform $\{a_1,a_2\}$ into
$\{b_1,b_2,b_3\}$ using, say, tree moves.  (We may write this more
briefly as ``one can get from $a_1,a_2$ to $b_1,b_2,b_3$ by tree
moves.'')  Then, for any multiset $A$, we can also transform
$A \cup \{a_1,a_2\}$ into $A \cup \{b_1,b_2,b_3\}$ using tree moves.
(Here $\cup$ is multiset union, where multiplicities are added.)
Such a transformation can then be used as part of further transformations.
Also, if one can get from $a_1,a_2$ to $b_1,b_2,b_3$ by tree
moves, then for any $c,d$ one can get from $r^c(1-r)^da_1,r^c(1-r)^da_2$
to $r^c(1-r)^db_1,r^c(1-r)^db_2,r^c(1-r)^db_3$ by tree moves (just multiply
every cylinder size involved by $r^c(1-r)^d$).  All of our
transformations of multisets will be sum-preserving.

Our next theorem shows that there are non-trivial examples of refinable 
numbers.
\begin{theorem}\label{selmerref}
If $r$ is the positive root of $x^n +x-1= 0$, where $n >1$ and $n
\not\equiv 5 \pmod 6$, then $r$ is refinable.
\end{theorem}

\begin{proof}
By a theorem of Selmer \cite{S}, the trinomial $x^n+x-1$ is
irreducible when $n
\not\equiv 5 \pmod 6$. Hence, $r$ is algebraic of degree $n$.
Let $A$ and $B$ be multisets of cylinder sizes for $r$ with the same
sum. Since $1-r = r^n$, we have
  $$
r^a(1-r)^{b+1} = r^{a+n}(1-r)^b.$$
The replacement of $r^a(1-r)^{b+1}$ by  $r^{a+n}(1-r)^b$ can be
thought of as both a trivial merge and a trivial split. Therefore, it
can be applied repeatedly to both $A$ and $B$ to produce new multisets
$A''$ and $B''$ containing only powers $r^i$, $i \geq 0$.

Next note that the replacement in a multiset
\begin{equation}\label{refine1}
r^a \rightarrow r^{a+1}, r^{a+n}
\end{equation}
is a split which is obtainable by tree moves (split $r^a$ to $r^{a+1}$
and $r^a(1-r)$ and merge $r^a(1-r)$ to $r^{a+n}$), so it can be
applied on both sides. Let $k$ be the largest exponent such that
$r^k$ occurs in $A''$ or $B''$. By repeatedly applying the
replacement (\ref{refine1})
to each $r^a$ with $a \leq k-n$, we can get from $A''$ and
$B''$ to $A'$ and $B'$ consisting entirely of powers $r^a$ with
$k-n+1 \leq a \leq k$. These $n$ powers of $r$ are linearly
independent over the rationals because $1,r,\ldots,r^{n-1}$ are
(since $r$ is algebraic of degree $n$).  So
the only way for $A'$ and $B'$ to have the same sum is to have
$A'$ = $B'$. This completes the proof.
\end{proof}

If $r$ is the positive root of $x^n+x-1 = 0$, and $s = r^d$ where
$d$ is a divisor of $n$, then $r \approx s$,
because $r = 1 - s^{n/d}$ (recall that, if $t$ is binomially reducible
to $r$, then so is $1-t$).  For most $n$, the preceding theorem shows
that $r$ is refinable.  Hence, to show that $r \sim_{top} s$,
it will suffice to show that $s$ is refinable.

In our next theorem we prove this for the special case when $n = 4$.

\begin{theorem}
If $r,s \in (0,1)$, $s = r^2$, and $r = 1-s^2$, then $r$ and $s$ are
refinable and the measures $\mu(r)$ and $\mu(s)$ are topologically equivalent.
\end{theorem}

\begin{proof}
We have $r = 1-r^4$, so $r \approx s$ (as noted above),
$r$ is refinable (by Theorem \ref{selmerref}), and
$r$ and $s$ are algebraic of degree $4$ (by Selmer's theorem
mentioned previously and the fact that binomial equivalence
preserves degree). Also, we have
\begin{equation}\label{eqaprime}
s = (1-s^2)^2 = (1-s)^2(1+2s+s^2).
\end{equation}
Now, from the cylinder size $1$, one can get to
\begin{equation}\label{eqa}
s^2,(1-s)^2,s^2(1-s),2s(1-s)^2,s^3(1-s),s^2(1-s)^2
\end{equation}
by tree splits (this corresponds to partitioning $\Cantor$ into
$\langle 1,1 \rangle$, $\langle 0,0 \rangle$,
$\langle 1,0,1 \rangle$, $\langle 1,0,0 \rangle$,
$\langle 0,1,0 \rangle$, $\langle 0,1,1,1 \rangle$,
and $\langle 0,1,1,0 \rangle$), and then to
\begin{equation}\label{eqb}
s,s^2,s^2(1-s),s^3(1-s)
\end{equation}
by a merge of the second, fourth, and sixth terms in (\ref{eqa}),
using (\ref{eqaprime}); of course,
one can also get from $1$ to
(\ref{eqb})
by an arbitrary split.  As noted earlier, this implies that
one can get from $s^a(1-s)^b$ to
$$
s^{a+1}(1-s)^b,s^{a+2}(1-s)^b,s^{a+2}(1-s)^{b+1},s^{a+3}(1-s)^{b+1}
$$
by an arbitrary split or by tree moves.

Next, one can get from $s$ to
\begin{eqnarray}\label{eqcprime}
s^3,s(1-s)^2,2s^2(1-s)^2,2s^4(1-s),s^3(1-s)^2,
\qquad \qquad \nonumber \\ \hfill
s^3(1-s)^3,s^5(1-s)^2,s^4(1-s)^3
\end{eqnarray}
by tree splits --- multiply the steps from $1$ to (\ref{eqa}) by $s$,
and then use three successive tree splits to replace the
cylinder size $s^3(1-s)$ with
$$s^4(1-s),s^3(1-s)^3,s^5(1-s)^2,s^4(1-s)^3.$$
One can then get from (\ref{eqcprime}) to
\begin{align}\label{eqbprime}
s^3,2s(1-s)^2,s^2(1-s)^2,2s^4(1-s),s^5(1-s)^2
\end{align}
by a merge from $s(1-s)^2$ times (\ref{eqb}) to $s(1-s)^2$. But we have
$$
s = r^2 = (1-s^2)^2 = (1-s)^2 + 2s(1-s)^2 + s^2(1-s)^2,
$$
and $s$ is also equal to the sum of the terms in (\ref{eqbprime})
since all of our moves are sum-preserving,
so we must have
$$(1-s)^2 = s^3 + 2s^4(1-s) + s^5(1-s)^2
$$
and hence we can get from (\ref{eqbprime}) to 
\begin{align}\label{eqc}
(1-s)^2,2s(1-s)^2,s^2(1-s)^2,
\end{align}
 by a merge.
So, we can get from
$s$ to (\ref{eqc}) by tree moves as well as by a split.

Next, we can get from $1$ to $s,(1-s)$ by a tree split and then to
\begin{equation}\label{eqd}
s,s(1-s),s^2(1-s),s^2(1-s)^2,s^3(1-s)^2
\end{equation}
by a (\ref{eqb}) move.  (Here ``by a (\ref{eqb}) move'' is short for
``by replacing a cylinder size $t$ with $t$ times (\ref{eqb}), which
is a split and can also be accomplished by tree moves.''  In this case
$t = 1-s$.)
Finally, we can get from $1$ to $s,(1-s)$ by a tree
split, then to
$$
(1-s),(1-s)^2,2s(1-s)^2,s^2(1-s)^2
$$
by a (\ref{eqc}) move, then to
\begin{equation}\label{eqe}
(1-s),3s(1-s)^2,2s^2(1-s)^2,s^2(1-s)^3,s^3(1-s)^3
\end{equation}
by a (\ref{eqb}) move.

Now let $A$ and $B$ be multisets of cylinder sizes with the same sum.
First use (\ref{eqc}) moves repeatedly on both multisets
to get rid of all cylinder sizes $s^a(1-s)^b$ with
$a > b+1$. Then use (\ref{eqb}) moves to get rid of all $s^a(1-s)^b$ with
$a<b$. So only numbers $s^a(1-s)^b$ with $a =b$ or $a = b+1$ occur in
the new multisets. Let $k$ be the largest exponent $a$ which occurs.
If either multiset contains a cylinder size $s^b(1-s)^b$ such that
$b < k-1$, then we can use a (\ref{eqd}) move to replace it
with cylinder sizes with larger exponents; similarly, we can use
a (\ref{eqe}) move to replace a cylinder size $s^{b+1}(1-s)^b$
such that $b < k-2$.  (Both of these steps yield new cylinder
sizes $s^{a'}(1-s)^{b'}$ with $a'=b'$ or $a'=b'+1$.)
By
performing these steps as many times as possible,
one can change all cylinder sizes to one 
of the
following five cylinder sizes:
\begin{eqnarray*}
s^{k-1}(1-s)^{k-2},\, s^{k-1}(1-s)^{k-1},\, s^k(1-s)^{k-1},
\qquad \qquad \\ \hfill
s^k(1-s)^k, \text{ or }s^{k+1}(1-s)^k.
\end{eqnarray*}
Then one can use tree splits on the first two of these five and
(\ref{eqc}) moves on the resulting occurrences of
$s^k(1-s)^{k-2}$ to reduce everything
to the cylinder sizes
$$
s^k(1-s)^{k-1},\,s^{k-1}(1-s)^k,\, s^k(1-s)^k, \text{ or }s^{k+1}(1-s)^k.
$$
Let $A'$ and $B'$ be the final multisets using these cylinder sizes
only. These four sizes are linearly independent over the
rationals. (One can verify this directly by noting that from $s$, $(1-s)$,
$s(1-s)$, $s^2(1-s)$ one can get $1$, $s$, $s^2$, $s^3$ as linear 
combinations, or
one can just notice that our argument shows that the multisets
$\{1\}$, $\{s\}$, $\{s^2\}$, $\{s^3\}$ can all be reduced to these four forms
using the same $k$. We are using here that $s$ has algebraic degree $4$.)
Therefore, since $A'$ and $B'$ are both multisets of these numbers and
$\sum A' = \sum A =\sum B = \sum B'$, we must have $A' = B'$. So $s$
is refinable, as desired. Finally, using Theorem \ref{newthm} we see that 
$\mu(r)$
and $\mu(s)$ are homeomorphic.
\end{proof}


\begin{corollary}
If $0 < r < 1$ and $r = 1-r^4$, then there are at least $4$ product
measures topologically equivalent to $\mu(r)$.
\end{corollary}

Let us mention some other results which can be proven using the
techniques of this paper.
It is not hard to show that any transcendental number is
refinable, although this is not useful for proving topological
equivalence.
One can show that, if $r$ is the root of $r^3 +r^2 -1 = 0$ in $(0,1)$
and $s = r^2$, then $r$ and $s$ are both refinable. This gives
another proof of the theorem of Navarro-Bermudez and Oxtoby. However,
it is simpler to produce the homeomorphism as they did. We have
verified that, if $r$ is the positive root of $r^6+r=1$,
then $r$, $r^2$, and $r^3$ are all refinable.
Thus, there are at least
six numbers in $(0,1)$ topologically equivalent to $r$.
We have also verified that the positive numbers $r$ and $s$
given by $s = r^4$ and $r = 1 - s^2$ are refinable.

A number of problems remain open.  One of these is the problem
stated after Definition~\ref{bineqdef}, which can now
be restated as follows.

\begin{problem}[{\cite[Problem 1065]{Ma}}]
Are $\sim_{top}$ and $\approx$ the same equivalence
relation on $[0,1]$?
\end{problem}

As noted earier, the above problem
has been solved in the negative by Austin~\cite{Au}.

In connection with this problem, we note that there are relatively
simple examples of two
probability measures on the Cantor space
each of which is a continuous image of the other but there is no
homeomorphism taking one to the other.  One such
exanple is $\mu(1/2)$ and $\nu$, where $\nu$ is
obtained from $\mu(1/2)$ by multiplying the measure of
any subset of the left half of the Cantor space by $3/2$
and multiplying the measure of any subset of the right half
of the Cantor space by $1/2$.  (Equivalently, $\nu$ can be
described as the disjoint sum of $(3/4) \mu(1/2)$ and
$(1/4) \mu(1/2)$.) The reason for this is that the $3/4$ half has no
clopen subset of measure $1/2^n$, since all the numerators are
divisible by $3$, whereas every clopen subset with
positive $\mu(1/2)$ measure has a clopen subset with measure $1/2^n$
for large $n$.

\begin{problem}[{\cite[Problem 1067]{Ma}}]
Is there an infinite  $\sim_{top}$ equivalence class?
Is there an infinite $\approx$ equivalence class?
\end{problem}

\begin{problem}
Is every number in a non-trivial $\sim_{top}$ equivalence
class refinable? (Are the corresponding measures good in Akin's sense?)
\end{problem}

The corresponding question about nontrivial $\approx$ equivalence classes
has a negative answer, by Theorem~\ref{newthm} combined with Austin's result.

Regarding this problem, we note that the number $1/3$ is not
refinable. This is because the partition $1/3 + 1/3 + 1/3$ of $1$
cannot be refined to a tree partition of $1$ since in any such tree
partition only one number will have odd numerator. In fact, one can
show that no rational $r$ in $(0,1)$ other than $1/2$ is refinable.

A particular case of interest in the preceding problem is the
remaining ``Selmer-like'' reals, and the numbers binomially
equivalent to Selmer or Selmer-like reals, where the Selmer reals
and the Selmer-like reals are the positive numbers $r$ satisfying
an equation $x^n+x-1=0$ where $n \not\equiv 5 \pmod 6$
or $n \equiv 5 \pmod 6$, respectively.
We already saw that the Selmer reals are refinable;
It turns out that one can show that the Selmer-like
reals are refinable as well, and the corresponding measures $\mu(r)$
are good in both cases.  But it remains open whether the numbers binomially
equivalent to them are refinable (and hence topologically equivalent
to them).

{\bf Acknowledgement} The authors would like to thank Mike Keane for
  his useful comments and information.


\begin{thebibliography}{10}

\bibitem{Ak} {\sc E. Akin},
Good measures on Cantor space, Trans. Amer. Math. Soc., to appear.

\bibitem{AP} {\sc S. Alpern and V. S. Prasad},
Typical dynamics of volume preserving homeomorphisms,
Cambridge Tracts in Mathematics, 139, Cambridge University Press,
Cambridge, 2000.

\bibitem{Au} {\sc T. D. Austin}, A pair of non-homeomorphic
measures on the Cantor set, preprint.

\bibitem{H} {\sc K.J. Huang},
Algebraic numbers and topologically equivalent measures in the Cantor
set, Proc. Amer. Math. Soc. 96 (1986), 560--562.

\bibitem{Ma}{\sc R. D. Mauldin}, Problems in topology arising from
   analysis, in Open problems in topology (J. van Mill and G. M. Reed,
   eds.), North-Holland,
   Amsterdam, 1990, pp.\ 617--629.

\bibitem{NB} {\sc F. J. Navarro-Bermudez},
Topologically equivalent measures in the Cantor space,
Proc. Amer. Math. Soc. 77 (1979), 229--236.

\bibitem{ONB} {\sc F. J. Navarro-Bermudez and J. C. Oxtoby},
Four topologically equivalent measures in the Cantor space,
Proc. Amer. Math. Soc. 104 (1988), 859--860.

\bibitem{O} {\sc J. C. Oxtoby}, Homeomorphic measures in metric spaces,
Proc. Amer. Math. Soc. 24 (1970), 419--423.

\bibitem{OP} {\sc J. C. Oxtoby and V. S. Prasad}, Homeomorphic
   measures in the Hilbert cube, Pac. J. Math. 77 (1978), 483--497.

\bibitem{OU}{\sc J. C. Oxtoby and S. M. Ulam},
  Measure preserving homeomorphisms and metrical transitivity,
Ann. Math. 42 (1941), 874--920.

\bibitem{P} {\sc R. G. E. Pinch}, Binomial equivalence of algebraic
   numbers, J. Indian Math. Soc. (N.S.) 58 (1992), 33--37.

\bibitem{S} {\sc E. S. Selmer}, On the irreducibility of certain
trinomials, Math. Scand. 4 (1956), 287--302.


\end{thebibliography}
\end{document}